\documentclass[reqno]{amsproc}%
\usepackage{amssymb}
\usepackage{amsmath}
\usepackage{eurosym}
\usepackage{amsfonts}
\usepackage{hyperref}%
\usepackage{graphicx}
\providecommand{\U}[1]{\protect\rule{.1in}{.1in}}
\advance \topmargin by -\headheight \advance
\topmargin by -\headsep
\newtheorem{theorem}{Theorem}
\newtheorem{definition}[theorem]{Definition}
\newtheorem{corollary}[theorem]{Corollary}
\newtheorem{lemma}[theorem]{Lemma}
\newtheorem{remark}[theorem]{Remark}
\numberwithin{equation}{section}
\begin{document}
\title[Generalized q-Morgan Voyce Polynomials]{Generalized q-Morgan Voyce Polynomials}
\author[Bahad\i r Y\i lmaz and Y\"{u}ksel Soykan]{Bahad\i r Y\i lmaz and Y\"{u}ksel Soykan}
\maketitle

\begin{center}
\textsl{Department of Mathematics, Faculty of Science, }

\textsl{Zonguldak B\"{u}lent Ecevit University, 67100, Zonguldak, Turkey}

\textsl{e-mail: \ yilmazbahadir.67@gmail.com}

\textsl{e-mail: \ yuksel\_soykan@hotmail.com}
\end{center}

\textbf{Abstract. }This study introduces and investigates generalized
q-Morgan-Voyce polynomials and their specific cases, including the first and
second kinds of q-Morgan-Voyce polynomials, q-Horadam-Morgan-Voyce
polynomials, and Fibonacci-type q-Morgan-Voyce polynomials. The generalized
q-Morgan-Voyce polynomials are defined by the recurrence relation featuring
the specific q-power $q^{(n-2)}$ and a negative sign, formulated as
$M_{n}(x,q)=(x+1+q)M_{n-1}(x,q)-q^{n-2}M_{n-2}(x,q)$ for $n\geq2$. The
generating functions and explicit expressions for these polynomials are
established by utilizing the Fibonacci operator and the binomial theorem.
Furthermore, the study extends these polynomials to negative indices and
derives the corresponding explicit formulas. Summation formulas and
determinantal presentations of the generalized polynomials and their special
cases are also comprehensively provided. Finally, the q-Cassini's formula for
the generalized q-Morgan-Voyce polynomials is established through the
definition of specific square matrices and their subsequent recurrence relations.

\textbf{2020 Mathematics Subject Classification \qquad}11B37,
11B50,\textbf{\ }11B39, 42C05, 33C45, 11C08.


\section{Introduction}

The polynomials $B_{n}(x)$ and $b_{n}(x)$, known as the first and second kind
Morgan-Voyce polynomials, respectively, were originally introduced by
Morgan-Voyce [\ref{MorganVoyce}] through the following recurrence relations:
\begin{align}
b_{n}(x) &  =x\,B_{n-1}(x)+b_{n-1}(x),\qquad(n\geq1),\label{eq:bn}\\
B_{n}(x) &  =(x+1)\,B_{n-1}(x)+b_{n-1}(x),\qquad(n\geq1),\label{eq:Bn}%
\end{align}
with the initial condition $b_{0}(x)=B_{0}(x)=1.$

These foundational definitions naturally give rise to the following
second-order linear recurrence relations:
\begin{equation}
b_{n}(x)=(x+2)b_{n-1}(x)-b_{n-2}(x)
\end{equation}
with $b_{0}(x)=1,$ $b_{1}(x)=x+1$ and
\begin{equation}
B_{n}(x)=(x+2)B_{n-1}(x)-B_{n-2}(x)
\end{equation}
with $B_{0}(x)=1,$ $B_{1}(x)=x+2.$

Swamy [\ref{Swamy1966}] extensively investigated the algebraic properties of
$B_{n}(x)$ and $b_{n}(x)$, subsequently applying these findings to derive
several noteworthy Fibonacci identities [\ref{Swamy1966fib}].

Later, Richard Andre-Jeannin [\ref{AndreJeannin1994}] introduced a generalized
family of Morgan-Voyce polynomials, denoted by $\{P_{n}^{(r)}\}$, which
satisfies the following second-order recurrence relation:
\begin{equation}
P_{n}^{(r)}(x)=(x+2)P_{n-1}^{(r)}(x)-P_{n-2}^{(r)}(x),\qquad n\geq
2,\label{A.Jreccurrance}%
\end{equation}
with initial conditions $P_{0}^{(r)}(x)=1\quad$and$\quad P_{1}^{(r)}%
(x)=x+r+1.$

Throughout this work, $r$ designates a fixed real parameter. From this
generalization, it is immediately apparent that
\[
P_{n}^{(0)}=b_{n},
\]
and that
\[
P_{n}^{(1)}=B_{n},
\]
which perfectly recover the classical Morgan-Voyce polynomials.

Similarly, A. F. Horadam [\ref{Horadam1995}] proposed another generalization,
denoted by $\{Q_{n}^{(r)}(x)\}$, governed by the recurrence:
\begin{equation}
Q_{n}^{(r)}(x)=(x+2)Q_{n-1}^{(r)}(x)-Q_{n-2}^{(r)}(x),\qquad(n\geq
2),\label{horadamrecurrence}%
\end{equation}
with initial conditions $Q_{0}^{(r)}(x)=2,$ $Q_{1}^{(r)}(x)=x+r+2.$

Furthermore, Horadam introduced two distinguished subclasses, now recognized
as the first and second kind Horadam-Morgan-Voyce polynomials. These are
characterized by the following respective second-order recurrence relations:
\[
C_{n}(x)=(x+2)C_{n-1}(x)-C_{n-2}(x)
\]
with $C_{0}(x)=2,$ $C_{1}(x)=x+2$ representing a Lucas-type polynomial
sequence, and
\[
c_{n}(x)=(x+2)c_{n-1}(x)-c_{n-2}(x)
\]
with $c_{0}(x)=1,$ $c_{1}(x)=x+3.$ It is worth noting the special parametric
connections: $P_{n-1}^{(2)}(x)=c_{n}(x)$ and $Q_{n}^{(0)}(x)=C_{n}(x)$.

Further properties of $C_{n}(x)$, particularly concerning their connections to
Lucas polynomials, are thoroughly detailed by Andre-Jeannin
[\ref{AndreJeannin1995}].

In a subsequent advancement, Swamy [\ref{Swamy2000}] formulated a broader
generalization of the modified Morgan-Voyce polynomials, denoted by
$\{w_{n}(a,b;x)\}=$ $\{w_{n}(x)\}$. This family, which may be referred to as
$p$-Morgan-Voyce polynomials, is constructed upon the following recurrence relation:

\begin{definition}
[\ref{Swamy2000}]For $n\geq2,$
\begin{equation}
w_{n}(x)=(x+p)w_{n-1}(x)-w_{n-2}(x),\qquad(n\geq
2),\label{reccurmodifiedmorgan}%
\end{equation}
with the initial conditions $w_{0}(x)=a\quad$and$\quad w_{1}(x)=b.$
\end{definition}

The sequence $w_{n}(x)$ can be naturally extended to negative indices via the
relation:
\begin{equation}
w_{-n}(x)=(x+p)w_{-\left(  n-1\right)  }(x)-w_{-\left(  n-2\right)  }(x)
\end{equation}

Consequently, the recurrence (\ref{reccurmodifiedmorgan}) remains valid for
all integers $n$. However, it is crucial to emphasize that for negative values
of $n$, this sequence does not strictly adhere to the standard definition of a polynomial.

Utilizing the recurrence relation (\ref{reccurmodifiedmorgan}), Horadam
[\ref{Horadam1965}] established the following structural identity:
\begin{equation}
w_{n}(x)=w_{1}(x)\widehat{U}_{n}(x)-w_{0}(x)\widehat{U}_{n-1}(x),
\end{equation}
where
\begin{equation}
\widehat{U}_{n}(x)=w_{n}(0,1;x).
\end{equation}

Before proceeding to our primary findings, we outline some fundamental
properties of the generalized modified Morgan-Voyce polynomials, an
overarching framework that encapsulates all classical variants. The
foundational recurrence relation yields the following second-order
characteristic equation:
\[
z^{2}-(x+p)z+1=0
\]
where the roots of this equation are $\alpha=\frac{x+p+\sqrt{(x+p)^{2}-4}}{2}$
and $\beta=\frac{x+p-\sqrt{(x+p)^{2}-4}}{2}.$ Throughout the remainder of this
study, we operate under the assumption that $\left\vert x+p\right\vert \neq2$.

Note that: For the sake of brevity, we shall refer to the generalizations of
the modified Morgan-Voyce polynomials as GM-Morgan-Voyce polynomials in the
subsequent sections.

Swamy [\ref{Swamy2000}] systematically identified five special cases of the
polynomials $w_{n}(x)$. These are designated as the first kind modified
Morgan-Voyce polynomial $\widetilde{b}_{n}(x)$, the second kind modified
Morgan-Voyce polynomial $\widetilde{B}_{n}(x)$, the first kind modified
Horadam-Morgan-Voyce polynomial $\widetilde{C}_{n}(x)$, the second kind
modified Horadam-Morgan-Voyce polynomial $\widetilde{c}_{n}(x)$, and
$\widehat{U}_{n}(x)$, which we shall refer to as the modified
Fibonacci-Morgan-Voyce polynomial. Their formal definitions are as follows:

\begin{definition}
[\ref{Swamy2000}]For $n\geq2,$ we have:
\[
\widetilde{b}_{n}(x)=(x+p)\widetilde{b}_{n-1}(x)-\widetilde{b}_{n-1}(x)
\]
with the initial conditions $\widetilde{b}_{0}(x)=1,$ $\widetilde{b}%
_{1}(x)=x+p-1,$
\[
\widetilde{B}_{n}(x)=(x+p)\widetilde{B}_{n-1}(x)-\widetilde{B}_{n-2}(x)
\]
with the initial conditions $\widetilde{B}_{0}(x)=1,$ $\widetilde{B}%
_{1}(x)=x+p,$
\[
\widetilde{C}_{n}(x)=(x+p)\widetilde{C}_{n-1}(x)-\widetilde{C}_{n-2}(x)
\]
with the initial conditions $\widetilde{C}_{0}(x)=2,$ $\widetilde{C}%
_{1}(x)=x+p,$
\[
\widetilde{c}_{n}(x)=(x+p)\widetilde{c}_{n-1}(x)-\widetilde{c}_{n-2}(x)
\]
with the initial conditions $\widetilde{c}_{0}(x)=1,$ $\widetilde{c}%
_{1}(x)=x+p+1,$
\[
\widehat{U}_{n}(x)=(x+p)\widehat{U}_{n-1}(x)-\widehat{U}_{n-2}(x)
\]
with the initial conditions $\widehat{U}_{0}(x)=0,$ $\widehat{U}_{1}(x)=1,$ respectively.
\end{definition}

We note that the first kind modified Horadam-Morgan-Voyce polynomial, denoted
by $\widetilde{C}_{n}(x)$, inherently exhibits a Lucas-type polynomial structure.

\subsection{Foundational Framework of $q$-Calculus}

To establish a rigorous mathematical background for the generalization of
polynomial sequences, it is essential to first outline the fundamental
principles of quantum calculus as systematically developed by Kac and Cheung
[\ref{Kac2002}]. Unlike classical analysis, quantum calculus eliminates the
necessity of taking limits, focusing instead on algebraic differences
parameterized by a quantum parameter $q \neq1$ [\ref{Kac2002}].

The foundational building block of this calculus is the $q$-differential
operator, which acts on an arbitrary function $f(x)$ as follows [\ref{Kac2002}%
]:
\begin{equation}
d_{q} f(x) = f(qx) - f(x)
\end{equation}
where $d_{q} x = (q-1)x$ [\ref{Kac2002}]. Consequently, the Jackson
$q$-derivative of a function $f(x)$ is defined as the plain ratio of these
quantum differentials [\ref{Kac2002}]:
\begin{equation}
D_{q} f(x) = \frac{d_{q} f(x)}{d_{q} x} = \frac{f(qx) - f(x)}{(q-1)x}, \quad x
\neq0
\end{equation}
The $q$-differentiation of a product of two functions lacks classical symmetry
and obeys the modified Leibniz rules [\ref{Kac2002}]:
\begin{equation}
D_{q} (f(x)g(x)) = f(qx)D_{q} g(x) + g(x)D_{q} f(x) = f(x)D_{q} g(x) +
g(qx)D_{q} f(x)
\end{equation}
Applying the operator $D_{q}$ to the monomial $x^{n}$ yields $D_{q} x^{n} =
[n]_{q} x^{n-1}$ [\ref{Kac2002}], where the key algebraic entity $[n]_{q}$
denotes the $q$-integer defined by [\ref{Kac2002}]:
\begin{equation}
[n]_{q} = \frac{q^{n} - 1}{q - 1} = 1 + q + q^{2} + \dots+ q^{n-1}%
\end{equation}
which clearly converges to the classical integer $n$ as $q \to1$
[\ref{Kac2002}]. This framework naturally extends to the $q$-factorial,
defined as $[n]_{q}! = [n]_{q} [n-1]_{q} \dots[1]_{q}$ with $[0]_{q}! = 1$
[\ref{Kac2002}].

In order to formulate expansions for more sophisticated polynomial families,
the classical binomial $(x-a)^{n}$ is replaced by its $q$-analogue, defined as
[\ref{Kac2002}]:
\begin{equation}
(x-a)_{q}^{n} = \prod_{k=0}^{n-1} (x - q^{k} a) = (x-a)(x-qa)\dots(x-q^{n-1}a)
\end{equation}
which satisfies the fundamental differentiation identity $D_{q} (x-a)_{q}^{n}
= [n]_{q} (x-a)_{q}^{n-1}$ [\ref{Kac2002}]. For any polynomial $f(x)$ of
degree $N$, a $q$-analogue of Taylor's formula can be established around a
point $c$ via the expansion [\ref{Kac2002}]:
\begin{equation}
f(x) = \sum_{j=0}^{N} (D_{q}^{j} f)(c) \frac{(x-c)_{q}^{j}}{[j]_{q}!}%
\end{equation}

By utilizing the $q$-Taylor expansion around $x=0$, one derives Gauss's
binomial formula for the polynomial $(x+a)_{q}^{n}$ [\ref{Kac2002}]:
\begin{equation}
(x+a)_{q}^{n} = \sum_{j=0}^{n} \binom{n}{j}_{q} q^{j(j-1)/2} a^{j} x^{n-j}%
\end{equation}
where the $q$-binomial coefficients are defined analogously to their classical
counterparts [\ref{Kac2002}]:
\begin{equation}
\binom{n}{j}_{q} = \frac{[n]_{q}!}{[j]_{q}! [n-j]_{q}!}%
\end{equation}
These coefficients exhibit deep geometric and combinatorial interpretations,
representing the exact number of $j$-dimensional subspaces within an
$n$-dimensional vector space over a finite field $\mathbb{F}_{q}$
[\ref{Kac2002}]. Furthermore, if two non-commuting operators or variables
satisfy the specific quantum commutation relation $yx = qxy$, the
noncommutative binomial formula can be stated as [\ref{Kac2002}]:
\begin{equation}
(x+y)^{n} = \sum_{j=0}^{n} \binom{n}{j}_{q} x^{j} y^{n-j}%
\end{equation}

Extending these structures to formal power series allows the derivation of
Heine's binomial formula [\ref{Kac2002}]:
\begin{equation}
\frac{1}{(1-x)_{q}^{n}}=\sum_{j=0}^{\infty}\frac{[n]_{q}[n+1]_{q}\dots\lbrack
n+j-1]_{q}}{[j]_{q}!}x^{j}=\sum_{j=0}^{\infty}\frac{(1-q^{n})_{q}^{j}%
}{(1-q)_{q}^{j}}x^{j}%
\end{equation}
Finally, taking the infinite limit $n\rightarrow\infty$ under the convergence
condition $|q|<1$ leads to Euler's two celebrated identities, which serve as
the foundations for the two distinct $q$-exponential functions [\ref{Kac2002}%
]:
\begin{equation}
e_{q}^{x}=\sum_{j=0}^{\infty}\frac{x^{j}}{[j]_{q}!}=\frac{1}{(1-(1-q)x)_{q}%
^{\infty}}\quad\left(  D_{q}e_{q}^{x}=e_{q}^{x}\right)
\end{equation}
\begin{equation}
E_{q}^{x}=\sum_{j=0}^{\infty}q^{j(j-1)/2}\frac{x^{j}}{[j]_{q}!}=(1+(1-q)x)_{q}%
^{\infty}\quad\left(  D_{q}E_{q}^{x}=E_{q}^{qx}\right)
\end{equation}
These foundational $q$-operators and non-commutative algebraic settings
provide the necessary tools to rigorously investigate the structural and
recurrence properties of the generalized polynomial sequences proposed in this study.

The mathematical study of $q$-polynomials has been profoundly shaped by
diverse operator-based and combinatorial methodologies. Carlitz
[\ref{Carlitz1975}] established the foundational basis by defining
$q$-Fibonacci polynomials through the combinatorial interpretation of binary
sequences with restricted consecutive ones. Cigler [\ref{Cigler1979},
\ref{Cigler1981}, \ref{Cigler2003}, \ref{Cigler2009}] systematically advanced
this field by applying umbral calculus and operator methods to derive
$q$-analogs for Fibonacci, Lucas, and Laguerre polynomials, revealing their
intricate connections to Morse code sequences, Bernoulli numbers, and
$q$-Genocchi structures. The significance of these sequences in broader
partition theory was further demonstrated by Andrews [\ref{Andrews2004}], who
successfully linked Fibonacci numbers to the Rogers-Ramanujan identities via
operator-theoretic proofs.

Recent research has expanded these horizons significantly: Jia, Liu, and Wang
[\ref{Jia2007}] provided comprehensive $q$-analogs for generalized Fibonacci
and Lucas polynomials using the Fibonacci operator $\eta_{x}$ approach, while
Belbachir and Benmezai [\ref{Belbachir2014}] introduced alternative $q$-Lucas
formulations to address classical identity generalizations. Specialized
polynomial families have also received considerable attention, including the
$q$-Chebyshev polynomials developed by Cigler [\ref{Cigler2012},
\ref{Cigler2013}] and the discrete $q$-Hermite polynomials explored in
[\ref{Cigler2013b}]. The scope has further widened to include combinatorial
structures like $q$-Narayana polynomials [\ref{Cigler2012b}], and the
implementation of $(p,q)$-Fibonacci finite operators by Polatl\i
\ [\ref{Polatli2023}] to yield new tridiagonal determinantal representations.
Furthermore, recent studies have integrated these sequences into diverse
fields, ranging from arithmetic properties of $q$-Pell numbers [\ref{Pan2006}]
and bi-periodic tiling interpretations [\ref{Ramirez2016}] to bicomplex number
systems [\ref{Aydin2022}] and their role in geometric function theory
[\ref{Alsoboh2025}]. These collective works underscore that operator theory
provides a robust, unified framework for analyzing the recurrence and
structural properties of diverse polynomial families within the quantum
calculus domain.

\section{Main Result}

In this section, we introduce the generalized $q$-Morgan-Voyce polynomials,
which constitute the fundamental basis of our main results. Extending the
classical Morgan-Voyce sequences through the lens of $q$-calculus provides a
broader framework for examining their algebraic properties. The following
definitions establish the general recurrence relation for this new family and
detail the specific initial conditions required to construct the five special
cases of the generalized $q$-Horadam-Morgan-Voyce polynomials.

\begin{definition}
For $n\geq2$, the generalized $q-$Morgan-Voyce polynomials are defined by:
\begin{equation}
M_{n}(x,q)=(x+1+q)M_{n-1}(x,q)-q^{n-2}M_{n-2}(x,q)\label{recurqmorganvoyce}%
\end{equation}
with arbitrary initial conditions $M_{0}(x,q)$ and $M_{1}(x,q).$
\end{definition}

Building upon this general framework, the subsequent definition formally
isolates and categorizes the five notable special cases of the generalized $q
$-Morgan-Voyce polynomials.

\begin{definition}
For $n\geq2$, we have the following equalities establishing the sub-families:

\begin{description}
\item[(a)] The first kind $q-$Morgan-Voyce polynomial is defined by
\[
b_{n}(x,q)=(x+1+q)b_{n-1}(x,q)-q^{n-2}b_{n-2}(x,q)
\]
with the initial conditions $b_{0}(x,q)=q$ and $b_{1}(x,q)=x+q.$

\item[(b)] The second kind $q-$Morgan-Voyce polynomial is defined by
\[
B_{n}(x,q)=(x+1+q)B_{n-1}(x,q)-q^{n-2}B_{n-2}(x,q)
\]
with the initial conditions $B_{0}(x,q)=q$ and $B_{1}(x,q)=x+1+q.$

\item[(c)] The first kind $q-$Horadam-Morgan-Voyce polynomial is defined by
\[
C_{n}(x,q)=(x+1+q)C_{n-1}(x,q)-q^{n-2}C_{n-2}(x,q)
\]
with the initial conditions $C_{0}(x,q)=1+q$ and $C_{1}(x,q)=x+1+q. $

\item[(d)] The second kind $q-$Horadam-Morgan-Voyce polynomial is defined by
\[
c_{n}(x,q)=(x+1+q)c_{n-1}(x,q)-q^{n-2}c_{n-2}(x,q)
\]
with the initial conditions $c_{0}(x,q)=q$ and $c_{1}(x,q)=x+1+q+q^{2}.$

\item[(e)] The Fibonacci type $q-$Morgan-Voyce polynomial is defined by
\[
G_{n}(x,q)=(x+1+q)G_{n-1}(x,q)-q^{n-2}G_{n-2}(x,q)
\]
with the initial conditions $G_{0}(x,q)=0$ and $G_{1}(x,q)=1.$
\end{description}
\end{definition}

With the foundational definitions established, we now turn to the analytical
representation of these sequences. The following theorem derives the
generating function for the generalized $q$-Morgan-Voyce polynomials.

\begin{theorem}
\label{genqmorganvoycegen} Let $f(z)=\sum_{n=0}^{\infty} M_{n}(x,q)z^{n}$ be
the generating function of the generalized $q-$Morgan-Voyce polynomials. Then,
we have:
\begin{equation}
\sum_{n=0}^{\infty}M_{n}(x,q)z^{n}=\frac{M_{0}(x,q)+\left(  M_{1}%
(x,q)-(x+1+q)M_{0}(x,q)\right)  z}{(1-(x+1+q)z+z^{2}\eta_{z})}%
\end{equation}
where $\eta_{z}$ represents the Fibonacci operator, [\ref{Andrews2004}],
defined as $\eta_{z}f(z)=f(qz).$
\end{theorem}

Proof. To establish this generating function, we proceed as follows:
\begin{align*}
(1-(x+1+q)z+z^{2}\eta_{z})f(z)  & =f(z)-(x+1+q)zf(z)+z^{2}\eta_{z}f(z)\\
& =f(z)-(x+1+q)zf(z)+z^{2}f(qz)\\
& =\sum_{n=0}^{\infty}M_{n}(x,q)z^{n}-(x+1+q)\sum_{n=0}^{\infty}%
M_{n}(x,q)z^{n+1}+\sum_{n=0}^{\infty}M_{n}(x,q)q^{n}z^{n+2}\\
& =\sum_{n=0}^{\infty}M_{n}(x,q)z^{n}-(x+1+q)\sum_{n=1}^{\infty}%
M_{n-1}(x,q)z^{n}+\sum_{n=2}^{\infty}M_{n-2}(x,q)q^{n-2}z^{n}\\
& =M_{0}(x,q)+M_{1}(x,q)z-(x+1+q)M_{0}(x,q)z\\
& \quad+\sum_{n=2}^{\infty}\underbrace{\left(  M_{n}(x,q)-(x+1+q)M_{n-1}%
(x,q)+q^{n-2}M_{n-2}(x,q)\right) }_{\text{equal to zero}}z^{n}%
\end{align*}
which directly yields the desired result.

By systematically applying the initial conditions from Definition 3 to the
main generating function derived in Theorem \ref{genqmorganvoycegen}, we
immediately obtain the generating functions for each special case.

\begin{corollary}
The following generating function formulas hold true:

\begin{description}
\item[(a)] $\sum_{n=0}^{\infty}b_{n}(x,q)z^{n}=\frac{q-\left(  -x+qx+q^{2}%
\right)  z}{(1-(x+1+q)z+z^{2}\eta_{z})}.$

\item[(b)] $\sum_{n=0}^{\infty}B_{n}(x,q)z^{n}=\frac{q-\left(  q-1\right)
\left(  q+x+1\right)  z}{(1-(x+1+q)z+z^{2}\eta_{z})}.$

\item[(c)] $\sum_{n=0}^{\infty}C_{n}(x,q)z^{n}=\frac{(1+q)-q\left(
q+x+1\right)  z}{(1-(x+1+q)z+z^{2}\eta_{z})}.$

\item[(d)] $\sum_{n=0}^{\infty}c_{n}(x,q)z^{n}=\frac{q+\left(  1+x-qx\right)
z}{(1-(x+1+q)z+z^{2}\eta_{z})}.$

\item[(e)] $\sum_{n=0}^{\infty}G_{n}(x,q)z^{n}=\frac{z}{(1-(x+1+q)z+z^{2}%
\eta_{z})}.$
\end{description}
\end{corollary}

Next, we establish an alternative expression for the generating function,
initiated from an arbitrary index $n$, which provides a more adaptable
analytical tool for subsequent derivations.

\begin{theorem}
Let $g(z)=\sum_{k=n}^{\infty}M_{k}(x,q)z^{k}$ represent the shifted generating
function of the generalized $q-$Morgan-Voyce polynomials. Then:
\begin{equation}
\sum_{k=n}^{\infty}M_{k}(x,q)z^{k}=\frac{\left(  M_{n}(x,q)-q^{n-1}%
M_{n-1}(x,q)z\right)  z^{n}}{(1-(x+1+q)z+z^{2}\eta_{z})}.
\end{equation}

\end{theorem}

Proof. Operating similarly to the previous theorem, we manipulate the
summation operator:
\begin{align*}
(1-(x+1+q)z+z^{2}\eta_{z})g(z)  & =g(z)-(x+1+q)zg(z)+z^{2}\eta_{z}g(z)\\
& =g(z)-(x+1+q)zg(z)+z^{2}g(qz)\\
& =\sum_{k=n}^{\infty}M_{k}(x,q)z^{k}-(x+1+q)\sum_{k=n}^{\infty}%
M_{k}(x,q)z^{k+1}+\sum_{k=n}^{\infty}M_{k}(x,q)q^{k}z^{k+2}\\
& =\sum_{k=n}^{\infty}M_{k}(x,q)z^{k}-(x+1+q)\sum_{k=n+1}^{\infty}%
M_{k-1}(x,q)z^{k}+\sum_{k=n+2}^{\infty}M_{k-2}(x,q)q^{k-2}z^{k}\\
& =M_{n}(x,q)z^{n}+M_{n+1}(x,q)z^{n+1}-(x+1+q)M_{n}(x,q)z^{n+1}\\
& \quad+\sum_{k=n}^{\infty}\underbrace{\left(  M_{k}(x,q)-(x+1+q)M_{k-1}%
(x,q)+q^{k-2}M_{k-2}(x,q)\right) }_{\text{equal to zero}}z^{k}\\
& =\left(  M_{n}(x,q)+\left(  M_{n+1}(x,q)-(x+1+q)M_{n}(x,q)\right)  z\right)
z^{n}\\
& =\left(  M_{n}(x,q)-q^{n-1}M_{n-1}(x,q)z\right)  z^{n}.
\end{align*}

Specializing this shifted generating function to the five specific families
yields the following corollary.

\begin{corollary}
The following expressions hold true:

\begin{description}
\item[(a)] $\sum_{k=n}^{\infty}b_{k}(x,q)z^{k}=\frac{\left(  b_{n}%
(x,q)-q^{n-1}b_{n-1}(x,q)z\right)  z^{n}}{(1-(x+1+q)z+z^{2}\eta_{z})}.$

\item[(b)] $\sum_{k=n}^{\infty}B_{k}(x,q)z^{k}=\frac{\left(  B_{n}%
(x,q)-q^{n-1}B_{n-1}(x,q)z\right)  z^{n}}{(1-(x+1+q)z+z^{2}\eta_{z})}.$

\item[(c)] $\sum_{k=n}^{\infty}C_{k}(x,q)z^{k}=\frac{\left(  C_{n}%
(x,q)-q^{n-1}C_{n-1}(x,q)z\right)  z^{n}}{(1-(x+1+q)z+z^{2}\eta_{z})}.$

\item[(d)] $\sum_{k=n}^{\infty}c_{k}(x,q)z^{k}=\frac{\left(  c_{n}%
(x,q)-q^{n-1}c_{n-1}(x,q)z\right)  z^{n}}{(1-(x+1+q)z+z^{2}\eta_{z})}.$

\item[(e)] $\sum_{k=n}^{\infty}G_{k}(x,q)z^{k}=\frac{\left(  G_{n}%
(x,q)-q^{n-1}G_{n-1}(x,q)z\right)  z^{n}}{(1-(x+1+q)z+z^{2}\eta_{z})}.$
\end{description}
\end{corollary}

To facilitate our algebraic manipulations in the operator calculus context, we
briefly recall a generalized form of the binomial theorem.

\begin{theorem}
\label{binomqalgebra} [\ref{Kac2002}] If $yx=qxy$ where $q$ is a number
commuting with both $x$ and $y$, then we have:
\begin{equation}
(x+y)^{n}=\sum_{j=0}^{n} \binom{n}{j}_{q} x^{n-j}y^{j}.
\end{equation}

\end{theorem}

Furthermore, we derive the following lemma. It establishes two key operator
identities that will serve as the fundamental algebraic scaffolding for the
remainder of this work.

\begin{lemma}
\label{lemmebinomqoperator}Let $A$ and $B$ be operators defined as follows:
\begin{align*}
A  & =(x+1+q)z,\\
B  & =-z^{2}\eta_{z},
\end{align*}
Then, the following properties hold:

\begin{description}
\item[(a)] For $n\geq0,$
\[
(A+B)^{n}=\sum_{j=0}^{n} \binom{n}{j}_{q} A^{n-j}B^{j}.
\]

\item[(b)] For $n\geq0,$
\[
B^{n}=(-1)^{n}z^{2n}q^{n(n-1)}\eta_{z}^{n}.
\]

\end{description}
\end{lemma}

Proof. (a) Let $h(z)$ be an arbitrary function. We evaluate the composition:
\begin{align*}
BA(h(z)) &  =B(A(h(z)))\\
&  =B((x+1+q)zh(z))\\
&  =-z^{2}\eta_{z}(((x+1+q)zh(z)))\\
&  =-z^{2}(x+1+q)qz\eta_{z}h(z)\\
&  =q(x+1+q)z\left(  -z^{2}\eta_{z}h(z)\right) \\
&  =q(x+1+q)zB(h(z))\\
&  =qA(B(h(z)))\\
&  =qAB(h(z)).
\end{align*}
This demonstrates that $BA=qAB.$ By applying Theorem \ref{binomqalgebra}, the
result is immediately obtained.

(b) For the proof, we proceed by induction on $n.$ For $n=0,$ the identity
trivially follows. Assuming the identity holds for $n=k,$ we evaluate the case
for $n=k+1$ on an arbitrary function $h(z)$:
\begin{align*}
B^{k+1}(h(z))  & =B\left(  B^{k}h(z)\right)  =B\left(  ((-1)^{k}%
z^{2k}q^{k(k-1)}\eta_{z}^{k})h(z)\right) \\
& =-z^{2}\eta_{z}((-1)^{k}z^{2k}q^{k(k-1)}\eta_{z}^{k}h(z))\\
& =-z^{2}(-1)^{k}(qz)^{2k}q^{k(k-1)}\eta_{z}^{k+1}h(z)\\
& =(-1)^{k+1}z^{2k+2}q^{k(k+1)}\eta_{z}^{k+1}h(z).
\end{align*}
This implies that
\[
B^{k+1}=(-1)^{k+1}z^{2k+2}q^{k(k+1)}\eta_{z}^{k+1}%
\]
which completes the induction step and the proof.

Combining the principles of the $q$-commuting binomial theorem and the
operator properties established in Lemma \ref{lemmebinomqoperator}, we now
define and expand the Morgan-Voyce operator.

\begin{theorem}
[Morgan-Voyce operator]Let $\Omega=(x+1+q)z-z^{2}\eta_{z}$. For $n\geq0$, we
have the explicit expansion:
\[
\Omega^{n}=\sum_{j=0}^{n}(-1)^{j} \binom{n}{j}_{q} (x+1+q)^{n-j}%
q^{j(j-1)}z^{n+j}\eta_{z}^{j}%
\]

\end{theorem}

Applying the Morgan-Voyce operator $\Omega^{n}$ to the specific monomials
$z^{2}$ and $z$ directly yields the following explicit actions:

\begin{corollary}
\label{corolmorganvoyveoperator}For $n\geq0,$ we obtain:

\begin{description}
\item[(a)]
\[
\Omega^{n}z^{2}=z^{n+2}\sum_{j=0}^{n}(-z)^{j} \binom{n}{j}_{q} (x+1+q)^{n-j}%
q^{j(j+1)}.
\]

\item[(b)]
\[
\Omega^{n}z=z^{n+1}\sum_{j=0}^{n}(-z)^{j} \binom{n}{j}_{q} (x+1+q)^{n-j}%
q^{j^{2}}.
\]

\end{description}
\end{corollary}

Armed with the operational machinery developed above, we are now positioned to
state and prove an explicit combinatorial expression for the generalized
$q$-Morgan-Voyce polynomials.

\begin{theorem}
\label{explicitformulagenmorganqankly}For $n\geq2$, the generalized
$q$-Morgan-Voyce polynomials satisfy the following explicit formulation:
\[
M_{n}(x,q)=M_{0}(x,q)A_{n}(x,q)+M_{1}(x,q)E_{n}(x,q).
\]
where
\begin{align*}
A_{n}(x,q)  & =\sum_{j=0}^{\left\lfloor \frac{n-2}{2}\right\rfloor }(-1)^{j+1}
\binom{n-j-2}{j}_{q} (x+1+q)^{n-2j-2}q^{j(j+1)},\\
E_{n}(x,q)  & =\sum_{j=0}^{\left\lfloor \frac{n-1}{2}\right\rfloor }(-1)^{j}
\binom{n-j-1}{j}_{q} (x+1+q)^{n-2j-1}q^{j^{2}}.
\end{align*}

\end{theorem}

Proof. Using Theorem \ref{genqmorganvoycegen}, we express the summation
dynamically:
\begin{align}
\sum_{n=0}^{\infty}M_{n}(x,q)z^{n}  & =\frac{1}{(1-(x+1+q)z+z^{2}\eta_{z}%
)}\left(  M_{0}(x,q)+\left(  M_{1}(x,q)-(x+1+q)M_{0}(x,q)\right)  z\right)
\nonumber\\
& =\sum_{n=0}^{\infty}\Omega^{n}\left(  M_{0}(x,q)+M_{1}(x,q)z-(x+1+q)M_{0}%
(x,q)z\right) \nonumber\\
& =\sum_{n=0}^{\infty}\Omega^{n}M_{0}(x,q)+M_{1}(x,q)\sum_{n=0}^{\infty}%
\Omega^{n}z-(x+1+q)M_{0}(x,q)\sum_{n=0}^{\infty}\Omega^{n}%
z.\label{expilicitpositivetheogen}%
\end{align}

Now, we isolate the $n=0$ term from the first summation to apply the operator
properly:
\begin{align*}
\sum_{n=0}^{\infty}\Omega^{n}M_{0}(x,q)  & = M_{0}(x,q) + \sum_{n=1}^{\infty
}\Omega^{n}M_{0}(x,q)\\
& = M_{0}(x,q) + \sum_{n=0}^{\infty}\Omega^{n}\left( \Omega M_{0}(x,q)\right)
\\
& = M_{0}(x,q) + \sum_{n=0}^{\infty}\Omega^{n}\left( (x+1+q)zM_{0}(x,q) -
z^{2}M_{0}(x,q)\right) \\
& = M_{0}(x,q) + (x+1+q)M_{0}(x,q)\sum_{n=0}^{\infty}\Omega^{n}z -
M_{0}(x,q)\sum_{n=0}^{\infty}\Omega^{n}z^{2}.
\end{align*}

Substituting this back into the generating function
(\ref{expilicitpositivetheogen}), the terms $(x+1+q)M_{0}(x,q)\sum
_{n=0}^{\infty}\Omega^{n}z$ systematically cancel each other, yielding:
\[
\sum_{n=0}^{\infty}M_{n}(x,q)z^{n}=M_{0}(x,q)+M_{1}(x,q)\sum_{n=0}^{\infty
}\Omega^{n}z-M_{0}(x,q)\sum_{n=0}^{\infty}\Omega^{n}z^{2}.
\]

Hence, applying Corollary \ref{corolmorganvoyveoperator} for the remaining two
sums provides the following expansions:

\begin{itemize}
\item For the $M_{1}(x,q)$ coefficient, taking $n+j+1=m$, i.e., $j\leq
n\Rightarrow j\leq\frac{m-1}{2}$, we obtain:
\begin{align*}
M_{1}(x,q)\sum_{n=0}^{\infty}\Omega^{n}z  & =M_{1}(x,q)\sum_{n=0}^{\infty}%
\sum_{j=0}^{n}(-1)^{j} \binom{n}{j}_{q} (x+1+q)^{n-j}q^{j^{2}}z^{n+j+1}\\
& =M_{1}(x,q)\sum_{m=1}^{\infty}\sum_{j=0}^{\left\lfloor \frac{m-1}%
{2}\right\rfloor }(-1)^{j} \binom{m-j-1}{j}_{q} (x+1+q)^{m-2j-1}q^{j^{2}}z^{m}%
\end{align*}

\item For the $M_{0}(x,q)$ coefficient, taking $n+j+2=m$, i.e., $j\leq
n\Rightarrow j\leq\frac{m-2}{2}$, we similarly find:
\begin{align*}
-M_{0}(x,q)\sum_{n=0}^{\infty}\Omega^{n}z^{2}  & =-M_{0}(x,q)\sum
_{n=0}^{\infty}\sum_{j=0}^{n}(-1)^{j} \binom{n}{j}_{q} (x+1+q)^{n-j}%
q^{j(j+1)}z^{n+j+2}\\
& =M_{0}(x,q)\sum_{m=2}^{\infty} \sum_{j=0}^{\left\lfloor \frac{m-2}%
{2}\right\rfloor }(-1)^{j+1} \binom{m-j-2}{j}_{q} (x+1+q)^{m-2j-2}%
q^{j(j+1)}z^{m}%
\end{align*}

\end{itemize}

Comparing the coefficients of $z^{m}$ for $m\geq2$ on both sides yields the
explicit expression for $M_{m}(x,q).$ Since $m$ is merely a dummy index,
replacing $m$ directly with $n$ provides the desired formula:
\[
M_{n}(x,q)=M_{0}(x,q)\sum_{j=0}^{\left\lfloor \frac{n-2}{2}\right\rfloor
}(-1)^{j+1} \binom{n-j-2}{j}_{q} (x+1+q)^{n-2j-2}q^{j(j+1)}+M_{1}%
(x,q)\sum_{j=0}^{\left\lfloor \frac{n-1}{2}\right\rfloor }(-1)^{j}
\binom{n-j-1}{j}_{q} (x+1+q)^{n-2j-1}q^{j^{2}}.
\]

Note that for $m=0$ and $m=1$, equating the coefficients of $z^{n}$ on both
sides trivially returns the initial conditions $M_{0}(x,q)$ and $M_{1}(x,q),$
respectively. This completes the proof.\bigskip

By substituting the specific initial conditions of our five target families
into Theorem \ref{explicitformulagenmorganqankly}, we explicitly construct
their respective formulas.

\begin{corollary}
For $n\geq2$, the following explicit formulas hold:

\begin{description}
\item[(a)] $b_{n}(x,q)=qA_{1}(x,q)+(x+q)A_{2}(x,q).$

\item[(b)] $B_{n}(x,q)=qA_{1}(x,q)+(x+1+q)A_{2}(x,q).$

\item[(c)] $C_{n}(x,q)=(1+q)A_{1}(x,q)+(x+1+q)A_{2}(x,q).$

\item[(d)] $c_{n}(x,q)=qA_{1}(x,q)+(x+1+q+q^{2})A_{2}(x,q).$

\item[(e)] $G_{n}(x,q)=A_{2}(x,q).$
\end{description}

where $A_{1}(x,q)$ and $A_{2}(x,q)$ are explicitly stated in Theorem
\ref{explicitformulagenmorganqankly}.
\end{corollary}

Furthermore, comparing these explicit expressions reveals profound algebraic
relationships between the different sub-families, summarized in the following corollary.

\begin{corollary}
For $n\geq2$, the following interrelations are verified:
\begin{align*}
B_{n}(x,q)-b_{n}(x,q)  & =G_{n}(x,q),\\
c_{n}(x,q)-B_{n}(x,q)  & =q^{2}G_{n}(x,q),\\
c_{n}(x,q)-b_{n}(x,q)  & =\left(  1+q^{2}\right)  G_{n}(x,q).
\end{align*}

\end{corollary}

Having established the properties for positive indices, we now extend the
generalized $q$-Morgan-Voyce polynomials to negative indices. The following
lemma dictates the backward recurrence relation.

\begin{lemma}
For $n\geq1,$ the sequence follows the backward recurrence relation:
\[
M_{-n}(x,q)=q^{n}(x+1+q)M_{-n+1}(x,q)-q^{n}M_{-n+2}(x,q)
\]

\end{lemma}

Building upon this backward recurrence, the subsequent theorem provides the
explicit combinatorial formula for the generalized $q$-Morgan-Voyce
polynomials evaluated at negative indices.

\begin{theorem}
For $n\geq1,$ we have the explicit representation:
\[
M_{-n}(x,q)=C_{n}(x,q)M_{0}(x,q)+D_{n}(x,q)M_{1}(x,q),
\]
where the coefficients are given by:
\begin{align*}
C_{n}(x,q)  & =\sum_{k=0}^{\left\lfloor \frac{n}{2}\right\rfloor }%
(-1)^{k}q^{\tbinom{n+1}{2}-k(n-k)} \binom{n-k}{k}_{q} (x+1+q)^{n-2k}\\
D_{n}(x,q)  & =\sum_{k=0}^{\left\lfloor \frac{n-1}{2}\right\rfloor }%
(-1)^{k+1}q^{\tbinom{n+1}{2}-k(n-k)} \binom{n-k-1}{k}_{q} (x+1+q)^{n-2k-1}%
\end{align*}

\end{theorem}

Proof. For the proof, we proceed using strong induction on $n.$ First, for the
base case $n=1$, we have:
\begin{align*}
M_{-n}(x,q)  & =C_{1}(x,q)M_{0}(x,q)-D_{1}(x,q)M_{1}(x,q)\\
& =q(x+1+q)M_{0}(x,q)-qM_{1}(x,q),
\end{align*}
which is evidently true. We now assume that the statement of the theorem holds
for all integers up to $m$ (i.e., $n\leq m$). Consequently, for the step
$n=m+1$, we obtain:
\begin{align*}
M_{-(m+1)}(x,q)  & =q^{m+1}(x+1+q)M_{-m}(x,q)-q^{m+1}M_{-\left(  m-1\right)
}(x,q)\\
& =q^{m+1}(x+1+q)\left(  C_{m}(x,q)M_{0}(x,q)-D_{m}(x,q)M_{1}(x,q)\right) \\
& \quad-q^{m+1}\left(  C_{m-1}(x,q)M_{0}(x,q)-D_{m-1}(x,q)M_{1}(x,q)\right) \\
& =\left(  (x+1+q)q^{m+1}C_{m}(x,q)-q^{m+1}C_{m-1}(x,q)\right)  M_{0}(x,q)\\
& \quad-\left(  (x+1+q)q^{m+1}D_{m}(x,q)-q^{m+1}D_{m-1}(x,q)\right)
M_{1}(x,q)
\end{align*}

Thus, expanding the individual components yields the following equalities:

\begin{itemize}
\item
\begin{align*}
(x+1+q)q^{m+1}C_{m}(x,q)  & =\sum_{k=0}^{\left\lfloor \frac{m}{2}\right\rfloor
}(-1)^{k}q^{\tbinom{m+1}{2}+m+1-k(m-k)} \binom{m-k}{k}_{q} (x+1+q)^{m-2k+1}\\
& =\sum_{k=0}^{\left\lfloor \frac{m+1}{2}\right\rfloor }(-1)^{k}%
q^{\tbinom{m+2}{2}-k(m-k+1)}q^{k} \binom{m-k}{k}_{q} (x+1+q)^{m-2k+1}.
\end{align*}

\item
\begin{align*}
q^{m+1}C_{m-1}(x,q)  & =\sum_{k=0}^{\left\lfloor \frac{m-1}{2}\right\rfloor
}(-1)^{k}q^{\tbinom{m}{2}-k(m-k-1)} \binom{m-k-1}{k}_{q} (x+1+q)^{m-2k-1}\\
& =\sum_{k=0}^{\left\lfloor \frac{m-1}{2}\right\rfloor }(-1)^{k}q^{\tbinom
{m}{2}+m+1-k(m-k-1)} \binom{m-k-1}{k}_{q} (x+1+q)^{m-2k-1}.
\end{align*}
By shifting the index $k\rightarrow k-1,$ we deduce:
\begin{align*}
q^{m+1}C_{m-1}(x,q)  & =\sum_{k=1}^{\left\lfloor \frac{m+1}{2}\right\rfloor
}(-1)^{k-1}q^{\tbinom{m}{2}+m+1-\left(  k-1\right)  (m-k)} \binom{m-k}%
{k-1}_{q} (x+1+q)^{m-2k+1}\\
& =\sum_{k=0}^{\left\lfloor \frac{m+1}{2}\right\rfloor }(-1)^{k-1}%
q^{\tbinom{m+2}{2}-k(m-k+1)} \binom{m-k}{k-1}_{q} (x+1+q)^{m-2k+1}.
\end{align*}

\item
\begin{align*}
(x+1+q)q^{m+1}D_{m}(x,q)  & =\sum_{k=0}^{\left\lfloor \frac{m-1}%
{2}\right\rfloor }(-1)^{k}q^{\tbinom{m+1}{2}+m+1-k(m-k)} \binom{m-k-1}{k}_{q}
(x+1+q)^{m-2k}\\
& =\sum_{k=0}^{\left\lfloor \frac{m}{2}\right\rfloor }(-1)^{k}q^{\tbinom
{m+2}{2}-k(m-k+1)}q^{k} \binom{m-k-1}{k}_{q} (x+1+q)^{m-2k}.
\end{align*}

\item
\[
q^{m+1}D_{m-1}(x,q)=\sum_{k=0}^{\left\lfloor \frac{m-1}{2}\right\rfloor
}(-1)^{k}q^{\tbinom{m}{2}+m+1-k(m-1-k)} \binom{m-k-2}{k}_{q} (x+1+q)^{m-2k-2}.
\]
By shifting the index $k\rightarrow k-1,$ we find:
\begin{align*}
q^{m+1}D_{m-1}(x,q)  & =\sum_{k=1}^{\left\lfloor \frac{m}{2}\right\rfloor
}(-1)^{k-1}q^{\tbinom{m}{2}+m+1-\left(  k-1\right)  (m-k)} \binom{m-k-1}%
{k-1}_{q} (x+1+q)^{m-2k}\\
& =\sum_{k=0}^{\left\lfloor \frac{m}{2}\right\rfloor }(-1)^{k-1}%
q^{\tbinom{m+2}{2}-k(m-k+1)} \binom{m-k-1}{k-1}_{q} (x+1+q)^{m-2k}.
\end{align*}

\end{itemize}

Combining these individual expressions simplifies the coefficients
considerably:
\begin{align*}
(x+1+q)q^{m+1}C_{m}(x,q)-q^{m+1}C_{m-1}(x,q)  & =\sum_{k=0}^{\left\lfloor
\frac{m+1}{2}\right\rfloor }(-1)^{k}q^{\tbinom{m+2}{2}-k(m-k+1)}\\
& \quad\times\left(  q^{k} \binom{m-k}{k}_{q} + \binom{m-k}{k-1}_{q} \right)
(x+1+q)^{m-2k+1}\\
& =\sum_{k=0}^{\left\lfloor \frac{m+1}{2}\right\rfloor }(-1)^{k}%
q^{\tbinom{m+2}{2}-k(m-k+1)} \binom{m-k}{k}_{q} (x+1+q)^{m-2k+1}\\
& =C_{m+1}(x,q).
\end{align*}

\begin{align*}
(x+1+q)q^{m+1}D_{m}(x,q)-q^{m+1}D_{m-1}(x,q)  & =\sum_{k=0}^{\left\lfloor
\frac{m}{2}\right\rfloor }(-1)^{k}q^{\tbinom{m+2}{2}-k(m-k+1)}\\
& \quad\times\left(  q^{k} \binom{m-k-1}{k}_{q} + \binom{m-k-1}{k-1}_{q}
\right)  (x+1+q)^{m-2k}\\
& =\sum_{k=0}^{\left\lfloor \frac{m}{2}\right\rfloor }(-1)^{k}q^{\tbinom
{m+2}{2}-k(m-k+1)} \binom{m-k}{k}_{q} (x+1+q)^{m-2k}\\
& =D_{m+1}(x,q).
\end{align*}
Therefore, we confirm the validity of the explicit representation for negative indices.

To further analyze the coefficient sequences introduced in the explicit
formulas, the following lemma establishes their individual recursive structures.

\begin{lemma}
For $n\geq2,$ the structural coefficients satisfy the following recurrences:

\begin{description}
\item[(a)]
\[
A_{n}(x,q)=(x+1+q)A_{n-1}(x,q)-q^{n-2}A_{n-2}(x,q)
\]
with the initial conditions $A_{0}(x,q)=1$ and $A_{1}(x,q)=0.$

\item[(b)]
\[
E_{n}(x,q)=(x+1+q)E_{n-1}(x,q)-q^{n-2}E_{n-2}(x,q)
\]
with the initial conditions $E_{0}(x,q)=0$ and $E_{1}(x,q)=1.$ Note that:
$E_{n}(x,q)=G_{n}(x,q).$

\item[(c)]
\[
C_{n}(x,q)=q^{n}(x+1+q)C_{n-1}(x,q)-q^{n}C_{n-2}(x,q)
\]
with the initial conditions $C_{0}(x,q)=1$ and $C_{1}(x,q)=q(x+1+q).$

\item[(d)]
\[
D_{n}(x,q)=q^{n}(x+1+q)D_{n-1}(x,q)-q^{n}D_{n-2}(x,q)
\]
with the initial conditions $D_{0}(x,q)=0$ and $D_{1}(x,q)=q.$
\end{description}
\end{lemma}

Proof. The proof of (c) and (d) are presented in the last theorem.
Additionally the proof of (a) and (b) can be done similarly.

Applying the negative index formula to the specific initial conditions yields
explicit representations for the five distinct special cases.

\begin{corollary}
For $n\geq1$, the following formulas dictate the negative indices of the
special families:

\begin{description}
\item[(a)] $b_{-n}(x,q)=qC_{n}(x,q)+\left(  x+q\right)  D_{n}(x,q).$

\item[(b)] $B_{-n}(x,q)=qC_{n}(x,q)+(x+1+q)D_{n}(x,q).$

\item[(c)] $C_{-n}(x,q)=(1+q)C_{n}(x,q)+(x+1+q)D_{n}(x,q).$

\item[(d)] $c_{-n}(x,q)=qC_{n}(x,q)+(x+1+q+q^{2})D_{n}(x,q).$

\item[(e)] $G_{-n}(x,q)=D_{n}(x,q).$
\end{description}
\end{corollary}

Comparing these specific negative-index formulas uncovers the reverse
algebraic connections among the families.

\begin{corollary}
For $n\geq1$, the following algebraic differences hold true at negative
indices:
\begin{align*}
B_{-n}(x,q)-b_{-n}(x,q)  & =G_{-n}(x,q),\\
c_{-n}(x,q)-B_{-n}(x,q)  & =q^{2}G_{-n}(x,q),\\
c_{-n}(x,q)-B_{-n}(x,q)  & =(1+q^{2})q^{2}G_{-n}(x,q).
\end{align*}

\end{corollary}

We now turn our attention to the summation identities. The following theorem
derives key finite sums involving the generalized $q$-Morgan-Voyce polynomials.

\begin{theorem}
For $n\geq2$, we have the following summation evaluations:

\begin{description}
\item[(a)]
\[
\sum_{k=0}^{n} (q^{k}-(x+q))M_{k}(x,q)=M_{1}(x,q)-(x+q)M_{0}(x,q)-M_{n+1}%
(x,q)+q^{n}M_{n}(x,q).
\]

\item[(b)]
\[
\sum_{k=1}^{n} (x+1+q)^{n-k}q^{k}M_{k}(x,q)=(x+1+q)^{n}M_{2}(x,q)-M_{n+2}%
(x,q).
\]

\end{description}
\end{theorem}

Proof. \ (a) Iterating the recurrence relation of the generalized
$q-$Morgan-Voyce polynomials yields the system of equations:
\begin{align*}
M_{2}(x,q)  & =(x+1+q)M_{1}(x,q)-q^{0}M_{0}(x,q),\\
M_{3}(x,q)  & =(x+1+q)M_{2}(x,q)-q^{1}M_{1}(x,q),\\
& \vdots\\
M_{n}(x,q)  & =(x+1+q)M_{n-1}(x,q)-q^{n-2}M_{n-2}(x,q).
\end{align*}
Summing these equations vertically, we isolate the cumulative structure:
\[
\sum_{k=2}^{n}M_{k}(x,q)=\sum_{k=1}^{n-1} (x+1+q)M_{k}(x,q)-\sum_{k=0}^{n-2}
q^{k}M_{k}(x,q).
\]
By strategically adding and subtracting boundary terms to align the summation
indices, we deduce:
\begin{align*}
\sum_{k=0}^{n}M_{k}(x,q)-M_{0}(x,q)-M_{1}(x,q)  & =\sum_{k=0}^{n}%
(x+1+q)M_{k}(x,q)-\sum_{k=0}^{n}q^{k}M_{k}(x,q)-(x+1+q)M_{n}(x,q)\\
& \quad-(x+1+q)M_{0}(x,q)+q^{n-1}M_{n-1}(x,q)+q^{n}M_{n}(x,q).
\end{align*}
Rearranging the terms produces the final sum formulation:
\[
\sum_{k=0}^{n} (q^{k}-(x+q))M_{k}(x,q)=M_{1}(x,q)-(x+q)M_{0}(x,q)-M_{n+1}%
(x,q)+q^{n}M_{n}(x,q).
\]

(b) Dividing the primary recurrence relation by ascending powers of $(x+1+q)$,
we express the relation as a telescoping difference for $k\geq0$:
\[
\frac{q^{k}}{(x+1+q)^{k}}M_{k}(x,q)=\frac{1}{(x+1+q)^{k-1}}M_{k+1}%
(x,q)-\frac{1}{(x+1+q)^{k}}M_{k+2}(x,q).
\]

Summing this telescoping series from $k=1$ to $n$, the intermediate terms
inherently collapse:
\begin{align*}
\sum_{k=1}^{n}\frac{q^{k}}{(x+1+q)^{k}}M_{k}(x,q)  & =\sum_{k=1}^{n}\left(
\frac{1}{(x+1+q)^{k-1}}M_{k+1}(x,q)-\frac{1}{(x+1+q)^{k}}M_{k+2}(x,q)\right)
\\
& =M_{2}(x,q)-\frac{1}{(x+1+q)^{n}}M_{n+2}(x,q).
\end{align*}
Multiplying both sides by the factor $(x+1+q)^{n}$ firmly establishes the
theorem's second claim:
\[
\sum_{k=1}^{n} (x+1+q)^{n-k}q^{k}M_{k}(x,q)=(x+1+q)^{n}M_{2}(x,q)-M_{n+2}%
(x,q).
\]

Transitioning from algebraic and summation properties, we next formulate a
determinantal representation for the generalized $q$-Morgan-Voyce polynomials,
expressing them as determinants of specific tridiagonal matrices.

\begin{theorem}
For $n\geq2,$ the sequence can be expressed determinantally as:
\[
M_{n}(x,q)=\det(\widetilde{D}_{n}(x,q))
\]
where the structured matrix $\widetilde{D}_{n}(x,q)$ is defined as:
\[
\widetilde{D}_{n}(x,q)=\left(
\begin{array}
[c]{cccccccc}%
M_{1}(x,q) & M_{0}(x,q) &  &  &  &  &  & \\
1 & x+1+q & q &  &  &  &  & \\
& 1 & x+1+q & q^{2} &  &  &  & \\
&  & 1 & x+1+q & \dots &  &  & \\
&  &  & 1 & \dots & \dots &  & \\
&  &  &  & \dots & \dots & \dots & \\
&  &  &  &  & \dots & \dots & q^{n-2}\\
&  &  &  &  &  & 1 & x+1+q
\end{array}
\right)  .
\]

\end{theorem}

Proof. For the proof, we proceed by strong induction on $n$. Evaluating the
determinant for the fundamental base case $n=2$ immediately confirms the
formula. Assuming the property holds for all dimensions up to $k$ ($n\leq k$),
we expand the determinant along the last row for the $(k+1)$-dimensional
matrix to obtain the recursive step:
\begin{align*}
\det(\widetilde{D}_{k+1}(x,q))  & =\left(  x+1+q\right)  \det\left(
\widetilde{D}_{k}(x,q)\right)  -q^{k-1}\det\left(  \widetilde{D}%
_{k-1}(x,q)\right) \\
& =\left(  x+1+q\right)  M_{k}(x,q)-q^{k-1}M_{k-1}(x,q)\\
& =M_{k+1}(x,q).
\end{align*}

By assigning the appropriate initial matrix entries to the first row, this
general determinantal form is readily specialized for our five specific
polynomial families.

\begin{corollary}
For $n\geq2,$ the determinants explicitly recover the specialized sequences:

\begin{description}
\item[(a)] $b_{n}(x,q)=\left\vert
\begin{array}
[c]{cccccccc}%
x+q & q &  &  &  &  &  & \\
1 & x+1+q & q &  &  &  &  & \\
& 1 & x+1+q & q^{2} &  &  &  & \\
&  & 1 & x+1+q & \dots &  &  & \\
&  &  & 1 & \dots & \dots &  & \\
&  &  &  & \dots & \dots & \dots & \\
&  &  &  &  & \dots & \dots & q^{n-2}\\
&  &  &  &  &  & 1 & x+1+q
\end{array}
\right\vert .$

\item[(b)] $B_{n}(x,q)=\left\vert
\begin{array}
[c]{cccccccc}%
x+1+q & q &  &  &  &  &  & \\
1 & x+1+q & q &  &  &  &  & \\
& 1 & x+1+q & q^{2} &  &  &  & \\
&  & 1 & x+1+q & \dots &  &  & \\
&  &  & 1 & \dots & \dots &  & \\
&  &  &  & \dots & \dots & \dots & \\
&  &  &  &  & \dots & \dots & q^{n-2}\\
&  &  &  &  &  & 1 & x+1+q
\end{array}
\right\vert .$

\item[(c)] $C_{n}(x,q)=\left\vert
\begin{array}
[c]{cccccccc}%
x+1+q & 1+q &  &  &  &  &  & \\
1 & x+1+q & q &  &  &  &  & \\
& 1 & x+1+q & q^{2} &  &  &  & \\
&  & 1 & x+1+q & \dots &  &  & \\
&  &  & 1 & \dots & \dots &  & \\
&  &  &  & \dots & \dots & \dots & \\
&  &  &  &  & \dots & \dots & q^{n-2}\\
&  &  &  &  &  & 1 & x+1+q
\end{array}
\right\vert .$

\item[(d)] $c_{n}(x,q)=\left\vert
\begin{array}
[c]{cccccccc}%
x+1+q+q^{2} & q &  &  &  &  &  & \\
1 & x+1+q & q &  &  &  &  & \\
& 1 & x+1+q & q^{2} &  &  &  & \\
&  & 1 & x+1+q & \dots &  &  & \\
&  &  & 1 & \dots & \dots &  & \\
&  &  &  & \dots & \dots & \dots & \\
&  &  &  &  & \dots & \dots & q^{n-2}\\
&  &  &  &  &  & 1 & x+1+q
\end{array}
\right\vert .$

\item[(e)] $G_{n}(x,q)=\left\vert
\begin{array}
[c]{cccccccc}%
1 & 0 &  &  &  &  &  & \\
1 & x+1+q & q &  &  &  &  & \\
& 1 & x+1+q & q^{2} &  &  &  & \\
&  & 1 & x+1+q & \dots &  &  & \\
&  &  & 1 & \dots & \dots &  & \\
&  &  &  & \dots & \dots & \dots & \\
&  &  &  &  & \dots & \dots & q^{n-2}\\
&  &  &  &  &  & 1 & x+1+q
\end{array}
\right\vert .$
\end{description}
\end{corollary}

To investigate the $q$-analogue of Cassini's identity, we introduce two
specific matrices, $Q_{n}(x,q)$ and $W_{n}(x,q)$, which will serve as the
primary algebraic tools.

\begin{definition}
For $n\geq1,$ we formulate the square matrices $Q_{n}(x,q)$ and $W_{n}(x,q)$
as follows:
\begin{align*}
Q_{n}(x,q)  & =\left(
\begin{array}
[c]{cc}%
x+1+q & -q^{n-1}\\
1 & 0
\end{array}
\right)  ,\\
W_{n}(x,q)  & =\left(
\begin{array}
[c]{cc}%
M_{n+1}(x,q) & -\sqrt{q}^{n-1}M_{n}(z,q)\\
M_{n}(x,q) & -\sqrt{q}^{n-2}M_{n-1}(z,q)
\end{array}
\right)
\end{align*}
where the transformed variable is defined as $z=\frac{x+1+q}{\sqrt{q}}-1-q.$
\end{definition}

The relationship between these matrices is strictly recursive, as formalized
in the following theorem.

\begin{theorem}
For $n\geq2,$ the matrix decomposition governs their product sequence:
\begin{equation}
W_{n}(x,q)=Q_{n}(x,q)W_{n-1}(x,q).
\end{equation}

\end{theorem}

Proof. By executing direct matrix multiplication, we systematically verify
that $W_{n}(x,q)=Q_{n}(x,q)W_{n-1}(x,q)$ for $n\geq2.$ Expanding the product
implies:
\begin{align*}
Q_{n}(x,q)W_{n-1}(x,q)  & =\left(
\begin{array}
[c]{cc}%
x+1+q & -q^{n-1}\\
1 & 0
\end{array}
\right)  \left(
\begin{array}
[c]{cc}%
M_{n}(x,q) & -\sqrt{q}^{n-2}M_{n-1}(z,q)\\
M_{n-1}(x,q) & -\sqrt{q}^{n-3}M_{n-2}(z,q)
\end{array}
\right) \\
& =\left(
\begin{array}
[c]{cc}%
M_{n+1}(x,q) & -\sqrt{q}^{n-2}\left(  x+1+q\right)  M_{n-1}(z,q)+q^{n-1}%
\sqrt{q}^{n-3}M_{n-2}(z,q)\\
M_{n}(x,q) & -\sqrt{q}^{n-2}M_{n-1}(z,q)
\end{array}
\right)  .
\end{align*}

Consequently, applying the fundamental recurrence relation
(\ref{recurqmorganvoyce}) to the upper-right entry, we deduce:
\begin{align*}
-\sqrt{q}^{n-2}\left(  x+1+q\right)  M_{n-1}(z,q)  & + q^{n-1}\sqrt{q}%
^{n-3}M_{n-2}(z,q)\\
& =-q^{n-1}\left(  \left(  \frac{x+1+q}{\sqrt{q}}\right)  M_{n-1}%
(z,q)-q^{n-2}M_{n-2}(z,q)\right) \\
& =-\sqrt{q}^{n-1}M_{n}(z,q).
\end{align*}
By substituting this simplification back into the matrix, we successfully
recover the identity:
\begin{equation}
Q_{n}(x,q)W_{n-1}(x,q)=W_{n}(x,q).
\end{equation}

\begin{remark}
By systematically iterating the relation $W_{n}(x,q)=Q_{n}(x,q)W_{n-1}(x,q)$
down to the base index $n=1,$ we generate the expansive product chain:
\begin{equation}
W_{n}(x,q)=Q_{n}(x,q)Q_{n-1}(x,q)Q_{n-2}(x,q)\dots Q_{1}(x,q)W_{1}%
(x,q)\label{cassiniidentitylemma}%
\end{equation}

\end{remark}

By capitalizing on the matrix product expansion detailed in the previous
remark, we now establish the $q$-Cassini's formula for the generalized
$q$-Morgan-Voyce polynomials.

\begin{theorem}
[$q-$Cassini's formula]For $n\geq2,$ the polynomials satisfy the invariant
difference identity:
\begin{align*}
&  \left(  M_{n+1}(x,q)M_{n-1}(z,q)-\sqrt{q}M_{n}(x,q)M_{n}(z,q)\right) \\
& =q^{\frac{(n-1)^{2}}{2}}\left(  M_{2}(x,q)M_{0}(z,q)-\sqrt{q}M_{1}%
(x,q)M_{1}(z,q)\right)
\end{align*}

\end{theorem}

Proof. By mapping the determinant operator across both sides of the matrix
decomposition (\ref{cassiniidentitylemma}), we establish the multiplicative
property:
\[
\det W_{n}(x,q)=\det Q_{n}(x,q)\det Q_{n-1}(x,q)\dots\det Q_{2}(x,q)\det
W_{1}(x,q).
\]
Since the determinant of the recursive multiplier evaluates to:
\[
\det Q_{n}(x,q)=q^{n-1}%
\]
we consolidate the powers of $q$ to yield:
\begin{align*}
\det W_{n}(x,q)  & =q^{n-1}q^{n-2}\dots q^{1}\det W_{1}(x,q)\\
& =q^{\frac{\left(  n-1\right)  n}{2}}\det W_{1}(x,q).
\end{align*}
Expanding the determinants of both matrices, this implies:
\begin{align*}
\det W_{n}(x,q)  & =-\sqrt{q}^{n-2}M_{n+1}(x,q)M_{n-1}(z,q)+\sqrt{q}%
^{n-1}M_{n}(x,q)M_{n}(z,q)\\
& =-\sqrt{q}^{n-2}\left(  M_{n+1}(x,q)M_{n-1}(z,q)-\sqrt{q}M_{n}%
(x,q)M_{n}(z,q)\right)
\end{align*}
and simultaneously:
\begin{align*}
q^{\frac{\left(  n-1\right)  n}{2}}\det W_{1}(x,q)  & =q^{\frac{\left(
n-1\right)  n}{2}}\left(  -\sqrt{q}^{-1}M_{2}(x,q)M_{0}(z,q)+M_{1}%
(x,q)M_{1}(z,q)\right) \\
& =-q^{\frac{\left(  n-1\right)  n-1}{2}}\left(  M_{2}(x,q)M_{0}(z,q)-\sqrt
{q}M_{1}(x,q)M_{1}(z,q)\right)
\end{align*}
Therefore, by equating these two explicit forms, we structurally balance the
equation to extract the required $q$-Cassini result:
\[
\left(  M_{n+1}(x,q)M_{n-1}(z,q)-\sqrt{q}M_{n}(x,q)M_{n}(z,q)\right)
=q^{\frac{(n-1)^{2}}{2}}\left(  M_{2}(x,q)M_{0}(z,q)-\sqrt{q}M_{1}%
(x,q)M_{1}(z,q)\right)  .
\]

Finally, applying the general $q$-Cassini's formula to the specific parametric
settings yields the corresponding Cassini identities for the five specialized families.

\begin{corollary}
For $n\geq2$, we have the following specialized $q$-Cassini identities:

\begin{description}
\item[(a)] $b_{n+1}(x,q)b_{n-1}(z,q)-\sqrt{q}b_{n}(x,q)b_{n}(z,q)=q^{\frac
{(n-1)^{2}}{2}}\left[  (x+q)\left(  (q-1)(x+1+q)+\sqrt{q}\right)
-q^{2}\right] . $

\item[(b)] $B_{n+1}(x,q)B_{n-1}(z,q)-\sqrt{q}B_{n}(x,q)B_{n}(z,q)=q^{\frac
{(n-1)^{2}}{2}}\left[  (q-1)(x+1+q)^{2}-q^{2}\right]  .$

\item[(c)] $C_{n+1}(x,q)C_{n-1}(z,q)-\sqrt{q}C_{n}(x,q)C_{n}(z,q)=q^{\frac
{(n-1)^{2}}{2}}\left[  q(x+1+q)^{2}-(1+q)^{2}\right]  .$

\item[(d)] $c_{n+1}(x,q)c_{n-1}(z,q)-\sqrt{q}c_{n}(x,q)c_{n}(z,q)=q^{\frac
{(n-1)^{2}}{2}}\left[  (x+1+q+q^{2})\left(  (q-1)(x+1+q)-q^{2}\sqrt{q}\right)
-q^{2}\right]  .$

\item[(e)] $G_{n+1}(x,q)G_{n-1}(z,q)-\sqrt{q}G_{n}(x,q)G_{n}(z,q)=-q^{\frac
{(n-1)^{2}}{2}}\sqrt{q}.$
\end{description}
\end{corollary}

\textbf{Funding Declaration}

The authors received no financial support for the research, authorship, and/or
publication of this article.

\end{document}